\documentclass{article}
\usepackage{latexsym,amssymb,amsmath}
\usepackage{amsfonts,lscape}

\begin{document}

\title{The History of Barbilian's Metrization Procedure}

\date{}

\author {Wladimir G. Boskoff, Bogdan D. Suceav\u{a}} 
\maketitle

\begin{abstract}
Barbilian spaces are metric spaces with a metric
induced by a
special procedure of metrization which is inspired by
the study of
the models of non-Euclidean geometry. In the present material we
discuss the
history of Barbilian spaces and the evolution of the
theory. We
point out that some of the current references to the
work done in
Barbilian spaces refer to Barbilian's contribution
from 1934,
while his construction has been largely extended in
four works
published in Romanian in 1959-1962.
\end{abstract}

{\bf Keywords:} Barbilian spaces\footnote{The present document is an updated version of the material presented on April 17, 2005, in the AMS  Special Session on History of Mathematics, AMS Sectional Conference \# 1007,
abstract \# 1007-01-39, Santa Barbara, California. A different material on the same theme, entitled {\it Barbilian Spaces: The History of a Geometric Idea}, will appear in {\it Historia Mathematica}. For terminology, see also Paul J. Kelly's classical article \cite{K1954}.}, axiomatic geometry,  non-Euclidean geometry,
hyperbolic geometry. 

\vspace{.5cm}

{\bf MSC Number} Primary 01A60 Secondary 54E35, 51M10, 51M99.

\vspace{.5cm}

 Several new developments in the area of complex analysis from the last decade shed new light 
on the work of
 Dan Barbilian (1895-1961). Our desire to present here the history of the idea of 
Barbilian metrization procedure
is motivated by the increasing attention that Barbilian's idea from 1934 receives in  recent references (see
e.g. \cite{Beardon1998,WB1993, WB1995, BGS2007,GH2000,H2003, H2003b, H2004, H2004b, H2005b,HI2005,HL2004, HPS2005,I2002, 
I2003a, I2003b, 
S1999}). 
We believe that a complete study of the mathematical
problem must have in view the contributions Barbilian published in
his last years, when, motivated by growing interest in his work,
he extended and generalized the theory. The authors of the present note have studied how Barbilian's metrization procedure can be used to naturally generate some important classes of metrics in \cite{BGS2007}. 

Any discussion of Barbilian's personality could start
by
mentioning that he is one of the major cultural
personalities of
the Romanian culture in the twentieth century, having
made
contributions not only in mathematics, but also in
literature. A. Cioranescu has no hesitation to write 
(\cite{C1981}, p.9]) that he is ``one of the greatest
Romanian poets of our century and perhaps the greatest of all." His
literary work, published under the pen name of Ion
Barbu (a name
borrowed from his paternal grandfather), received the
utmost
consideration in Romanian literature, where he is
regarded today as one of the influential masters in Romanian poetry. Most
of his literary work was written in the twenties
and at the
beginning of the thirties. His most important literary
volume,
{\it Joc secund (Secondary Game)}, was published in 1930 and inspired
many critical
comments and studies (to mention just a few, see \cite{C1966,C1981,N1960,V1935}).

A. Cioranescu writes (\cite{C1981}, pp.35): ``Ion Barbu's biography significantly 
interweaves the urge
toward mathematical thought with interest in poetry. The combination of the two 
imaginations in one mind
is especially remarkable and probably unique."

 Several of Barbilian's biographers agree that the
year 1933 was his year of ``return to mathematics". He
was
interested in Euclidean geometry since his early
years, and when
he ``came back" (as he describes later on) to
mathematics, he
dedicated his time to ideas inspired from plane
geometry. He was
interested in generalizing the metric used in
non-Euclidean
geometry, more precisely in the Beltrami-Klein model. Already one year later, in 1934, he had some results
worth
being presented in an international meeting, in Prague.
The talk from
Prague is published as a two-page note \cite{B1934}. This
was followed
by an exchange of letters with
 Wilhelm Blaschke, who attended the Prague meeting and
was interested in Barbilian's metrization idea. Surprisingly,
D. Barbilian did not publish any further research paper on this
topic in the next five years, although he pursued this idea
thoroughly, and we know that there existed a joint project with
Blaschke. He also made a presentation in a conference of {\it
Societatea de \c{S}tiin\c{t}e Matematice}, on December 17, 1934,
and he included the idea of the metric as logarithmic oscillation in his lecture notes
\cite{B1933-1934}, pp. 405-409. From the historical standpoint,
it is useful to look for the source of the idea in his textbook, since his lectures contained many 
research ideas. To endorse this claim, we can
cite here a first hand witness: N. Mih\u{a}ileanu. He wrote in
(\cite{B1933-1934}, p. 10), referring to the period starting in
1941, relevant for Barbilian's approach for his classes and in
general his work:

\begin{quotation}
 {\it As his first assistant, in the
fall of
1941, [...] Barbilian sent me sometimes in the morning
the written
version of his lecture, so that I would prepare my
recitation. In
the mean time, [...]  he changed his mind and at 3
p.m. he was
coming to present his lecture, in a form which was not similar
any more,
except in some ideas, with the material written in the
same day.
[...] Thus, every lesson was an original piece of work
[...].}
\end{quotation}

Therefore, it is not unusual that we find his research idea  published in \cite{B1934}
in the lecture notes of the course
of
descriptive geometry written in the same year. It can be found in the
chapter
called {\it The induction problem,} after the section
dedicated to
the presentation of {\it the direct definition of the
non-Euclidean angle and distance}. The section is entitled
{\it The
variative definition of the non-Eucludean distance.}
(``Variative"
is a term that Barbilian made up, as he made up many
other terms,
when he needed them to better describe some concepts.)

In any case, after the Prague meeting, in a volume by
Leonard
Mascot Blumenthal, which appeared in 1938 at the University
of Missouri
\cite{B1938}, a section is dedicated to Barbilian
spaces. Thus,
the terminology was introduced by Blumenthal,
based on
Barbilian's contribution \cite{B1934}.
Blumenthal defended his Ph.D. thesis, entitled {\it
Lagrange
Resolvents in Euclidean Geometry} in 1927, at John
Hopkins
University. His advisor was Frank Morley. Over the
years, among
Blumenthal's many Ph. D. students at the University of
Missouri Columbia,
were L. M. Kelly (who later had a long
career at
Michigan State University) and William Kirk (who
worked at the
University of Iowa). It is clear that L. M.
Blumenthal's works
were influential and garnered much interest at
the time. The 1938 book, entitled {\it Distance
geometries; a
study of the development of abstract metrics,}
appeared with a
foreword by Karl Menger. 

In Bucharest, after 1939, Barbilian's mathematical
interest
shifted towards algebra and number theory.
Unfortunately, this
evolution took place before he managed to write down
in a final
form the theory and the results from 1934-1939. After
1940, there
was a period of political challenge and turmoil
throughout Europe
which impacted all aspects of academic life. One of the
consequences of the cold war and the appearing of the iron
curtain was
the disruption of personal and institutional contacts
in academic
life. For many original ideas, as is the case with
Barbilian
spaces, this situation interfered with the evolution
of the
theory. 

After the war, the mathematician who reintroduced the
problem of
Barbilian spaces was Paul Joseph Kelly. Born in 1915
in Riverside,
California, P. J. Kelly defended his thesis in 1942 at the
University
of Wisconsin - Madison. Then he served in the U.S.A.F.
during the
war. Shortly after the war, he worked as Instructor at
University
of Southern California, then he joined the faculty at the
University
of California Santa Barbara, from where he retired in
1982.
After 1957 he was for a few years the Chair
of the
Department of Mathematics at U. C. Santa Barbara.
Perhaps one of
his best known works is the volume  {\it
Projective geometry
and projective metrics,} written with Herbert
Busemann. Among
other works he wrote in the fifties, there is the
paper {\it
Barbilian geometry and the Poincar\'{e} model,}
published in {\it the
American Mathematical Monthly,} in the May 1954 issue \cite{K1954}.
In this
excellent expository paper, Kelly wrote:

\begin{quotation}
{\it In a very brief paper in \v{C}asopis
Matematiky
a Fysiky (1934-1935), D. Barbilian defines a class of
metric
spaces (which we will call Barbilian spaces) and
states some of
their properties. It seems to the author that the
Barbilian
approach to the Poincar\'e model has certain
advantages of
simplicity and generality. }
\end{quotation}

Thus, the paper presented in Prague in 1934 was still
of interest
twenty years later. P. J. Kelly published several
times in the
{\it Monthly} in the fifties and we can trace his interest
in
 Barbilian spaces until 1981 \cite{KM1981}.

In Romania, Barbilian must have found out about P. J. Kelly's
paper from the {\it Monthly} and the news surely made a strong
impression on him. He thought his research work from before the
war had been forgotten. He decided to return to geometry and,
since the notes were lost, he had to redo the whole theory from
scratch. It is clear from the first paragraph of the fragment we
will present below that P. J. Kelly's paper represented for
Barbilian an important motivation, maybe the major motivation that
made him come back to geometry. He worked through the details and wrote them in a publishable form. The
first paper from this sequence of works is {\it Asupra unui
principiu de metrizare,} ({\it On a metrization principle}), and
it appeared in the Bucharest journal {\it Studii \c{s}i
Cercet\u{a}ri Matematice}. But here
is the story, as Barbilian tells it, in the introduction of
\cite{B1959a}:

\begin{quotation}
{\it Recently, P. J. Kelly, in his article
entitled {\bf
Barbilian Geometry and the Poincar\'{e} model},
 reopens the discussion on the topic of our talk presented to the Meeting of mathematicians from slavic countries from Prague (1934), that he learned from L. M. Blumenthal's monograph, from the section entitled Barbilian spaces. Thus, this research regards a work done twenty four years ago. Since then, it was out of our interest. But between 1934 and 1939 we have speculated a lot on this idea, without publishing anything. We have communicated part of our results to Wilhelm Blaschke with whom, after the Prague Meeting, where he participated, we have exchanged letters on this topic. The long article we have projected, was not written in the end, due to the war and due to our change of our interest toward algebra and number theory. Since, as it seems, the geometries that we have proposed in 1934 are of interest today, we believe that it is not pointless to present the supplementary results obtained in this direction in the interval 1934-1939, reconstructing them after the brief notes that we kept. We are able to see that the viewpoint of the 1934 generalization can be surpassed by far. We will find a very general procedure of metrization, through which 
the positive functions of two points, on certain sets, can be refined to a distance.}
\end{quotation}

However, it seems that {\it the Monthly} was not available to Barbilian at
that time. He mentioned that he
 writes the work without being able to see
Kelly's paper:

\begin{quotation}
{\it We are not able to read P. J. Kelly's article,
which is
known to us only from the review written by L. M.
Blumenthal in
{\it Mathematical Reviews}. }
\end{quotation}

Blumenthal's review \cite{Blum1954} consists in just a few lines, and
for
completeness we reproduce it here. It seems this is all
Barbilian had
available to comment on Kelly's work.
\vspace{.1cm}

\begin{quotation} {\it Barbilian space [\v Casopis P\v est. Mat. Fys.
64,
182--183 (1935)] is obtained by assigning to each two
points $a,b$
of the interior $K$ of a simple closed, plane curve
$J$ the number
$$ d(a,b)=\log\,[\max_{p\in
J}\,(pa/pb)]+\log\,[\max_{q\in J}\,(qb/qa)] $$ as
distance, where
$xy$ denotes the euclidean distance of points $x,y$.
[A brief
investigation of this space is given by the reviewer,
Univ.
Missouri Studies 13, no. 2 (1938).] The author
establishes some
known properties of Barbilian space, discusses
Barbilian
geodesics, and defines a generalization in which $K$
and $J$ are
replaced, respectively, by any planar set and a closed
set $J$
disjoint from it.}
\end{quotation}

Barbilian travelled in 1957 to the Dresden
International Congress
of Mathematics. He must have had some access to
information, but
it seems that the works that interested him the most
were hard, if
not impossible, to find. 

Why did Barbilian write his four papers from the
interval
1959-1961 in Romanian? Why didn't he write all his
material in either
English, French or German? He was definitely
proficient in several
languages. In the twenties, when he studied at
G\"ottingen and planned a Ph.D. thesis under Edmund Landau's direction, he
attended classes taught by David Hilbert and Emily
Noether, among
other famous professors. His spouse was German and
there could
have been no problem in writing all his work in either
one of the
languages mentioned above. Why didn't he send his
material to
another journal instead of  {\it
Studii \c{s}i cercet\u{a}ri matematice}? We can only
speculate
about this decision. One possible hypothesis is that
he wanted to
ensure that  the
 material made it into print during his lifetime and that a local
mathematical journal with worldwide circulation was
the most
practical solution in the academic and political
context of
Romania in the late fifties.
In any case, this decision is one of the main reasons
why, later
on, some
 information was missing from the development of the
theory on Barbilian's metrization procedure. At least one person was well informed, though.
It was L. M.
Blumenthal, who wrote the comment for {\it
Mathematical Reviews.}
We cite here the whole review from Mathematical
Reviews, as it
includes the description of Barbilian's metrization
principle \cite{Blum1959}:

\begin{quotation}

{\it In an earlier paper [\v{C}asopis P\v est. Mat.
Fys. 64
(1934/35), 182--183] the author introduced and stated
without
proof some properties of the metric space obtained by
attaching to
each two points $a,b$ of the interior of a simple
closed plane
curve $K$ the distance $$ d(a,b)=\max_{p,q\in
K}\ln[(pa/pb)(qb/qa)].
$$ The present paper greatly extends this metrization
procedure and investigates in detail the resulting
space. Let $K$
denote a closed subset of a topological space $T,J$
any subset of
$T$ and $(PA)$ a real positive function of points $P,A
(P\in
K,A\in J)$ which is continuous with respect to $P$ for
$A$ fixed.
Suppose the ratio $(PA)/(PB)$ is constant only when
$A=B$ ($P$ a
variable element of $K$ and $A,B$ fixed elements of
$J$). Putting
$M=\max_{P\in K}(PA)/(PB)$, $m=\min_{P\in
K}(PA)/(PB)$, define
$d(A,B)=\ln M/m (A,B\in J)$. The paper studies the
geometry of the
resulting metric space.}
\end{quotation}

While Barbilian
continued to produce papers in this direction, some of
the later
references still emphasized his 1934 work. He
wrote a second
paper entitled {\it Fundamentele metricilor abstracte
ale lui
Poincar\'e \c{s}i Carath\'eodory ca aplica\c{t}ie a
unui principiu
general de metrizare} ({\it Foundations of the
abstract metrics of
Poincar\'e and Carath\'eodory as application of a
general
principle of metrization} (a work presented at the
Institute of
Mathematics of Bucharest on the 4th of June 1959),
published also
in Romanian in the journal {\it Studii \c{s}i
cercet\u{a}ri
matematice,} in the same year \cite{B1959b}. The
next work
is {\it J-metricile naturale finsleriene} ({\it The
Finslerian
natural J-metrics}) \cite{B1960}, and the
work written
jointly with Nicolae Radu, which appeared posthumuously,
{\it
J-metricile naturale finsleriene \c{s}i func\c{t}ia de
reprezentare a lui Riemann (Finslerian natural J-metrics and Riemann's 
reprezentation 
function} \cite{BR1962}. In fact,
the last
work was submitted to the editor on October 20,
1961, after
Barbilian's death, which was on August 11, 1961.

There may be of interest to show here the most general
part of the
construction done by Barbilian in 1959 
\cite{B1959a}. In our
English translation, the construction is the
following.

\begin{quotation}
{\it
{\bf Definition.} {\it Let $K,J$ be two arbitrary sets
and $(PA)$
a function of the pair $P\in K,$ $A \in J,$ with real
and positive
values. We call $(PA)$ the influence of the set $K$
over the set
$J.$}
The only hypothesis satisfied by the influence is the
following:
{\it The extremum requirement.} For $A,B \in J$ fixed
and $P \in
K$ variable, the ratio of influences
\begin{equation}
\label{ratio} (PA) : (PB)
\end{equation}
 reaches a maximum $\mathcal{M}$ (obviously finite and
positive).
This yields the following.
{\it Consequence.} The ratio of influences
(\ref{ratio}) also
reaches a minimum $m.$
Indeed, due to the extremum requirement, the ratio of
influences
$(PB):(PA)$ reaches a positive and finite maximum,
that we may
denote $1/m.$ Clearly, therefore, that (\ref{ratio})
reaches a
minimum (finite and positive) $m.$
Here are a few constructions of the requirement of
extremum:

1. {\it Our construction from Prague.}
$K=$ a closed simple curve (Jordanian curve). $J=$ the
interior
domain of $K.$ The influence = the Euclidean distance
$(PA).$

2. {\it Dr. P. J. Kelly's construction.}
$K=$ a closed set in the Euclidean plane. $J=$ a
disjoint set,
i.e.
\begin{equation}
J \cap K = \emptyset
\end{equation}
but otherwise $J$ is arbitrary in the same Euclidean
plane. The
influence = the Euclidean distance $(PA).$

3. {\it Our construction from November 1934}
(Communicated by W.
Blaschke in a letter from the same period).
$K=$ a closed set in the Euclidean space with $n$
dimensions. $J=$
an abstract arbitrary set. The influence = a
pseudo-distance (that
is, a nonnegative symmetric function of the part
$X,Y$, vanishing
if $X=Y$) $(PA)$ defined on the abstract Cartesian
product $ K
\times J$ and continuous on $K,$ with respect to the
argument $P
\in K.$

4. {\it The new construction.}
$K=$ a compact set in a topological space $T$, $K
\subset T,$ $J=$
an arbitrary abstract set. The influence = a positive
function
$(PA),$ $P \in K, $ $A\in J,$ continuous on $K$ (with
respect to
$P$), but otherwise arbitrary.
}
\end{quotation}

Then, Barbilian states and proves the fundamental
theorem.
\vspace{.2in}

\begin{quotation} {\it  {\bf Fundamental Theorem.} Let (PA) be
the influence
of the set $K$ over an arbitrary set $J$, satisfying
the extremum
requirement. Then, the logarithmic oscillation
$$AB = \log \frac{\mathcal{M}}{m}$$
of the ratio of influences $\frac{(PA)}{ (PB)},$ with
$A,B \in J$
fixed and $P \in K$ variable, defines in $J$ a weak
distance.
}
\end{quotation}

For Barbilian, a weak distance $d$ satisfies the
symmetry
condition, the triangle inequality, and $A=B$ implies
$d(A,B) =0,$
but not the converse. This is the idea of the
construction
done in \cite{B1959a}.

 In 1981 \cite{KM1981}, P. J. Kelly
and G. Matthews pointed out how the idea of Barbilian space can be
integrated in the study of non-Euclidean geometry at an
introductory level. They have presented Kelly's idea from 1954 in an accesible style
suitable for a textbook. The study of the connection between the
Apollonian metric and quasiconformal mappings attracted the
interest of many researchers in recent years, especially since the
topic reappeared naturally in an area traditionally lying in the
field of analysis.  We mention here just a few of the papers that cite Barbilian's work.
A first study in this direction is due to A.
F. Beardon \cite{Beardon1998}.  F. W. Gehring and K. Hag, in
\cite{GH2000}, have continued the investigation. P. H\"ast\"o
\cite{H2003, H2003b, H2004, H2004b, H2005b}, Z. Ibragimov
\cite{I2002, I2003a, I2003b},   P. H\"ast\"o and Z. Ibragimov \cite{HI2005},  
P. A. H\"ast\"o and 
H. Lind\'en \cite{HL2004},
P. H\"ast\"o, 
S. Ponnusamy and S. K. Sahoo \cite{HPS2005} have detailed the study of
Apollonian metric.  A comprehensive monograph in the topic, from
the geometric viewpoint, including a discussion of Barbilian's
concepts from 1959-1960 is due to W. G. Boskoff
\cite{WB1996carte}, and recent developments in the geometry of various classes of metrics
are due to W. G. Boskoff, M. G. Ciuc\u{a} and B. D. Suceav\u{a} \cite{BGS2007}.  All these works cite and have a common
source in Barbilian's paper presented in Prague, and they convey
the potential lying in this research direction.

Regarding the term ``Apollonian metric (space),"
Barbilian himself
 suggested in 1959 that this term would be best,
pointing out
with modesty that he finds it more appropriate than
the name
``Barbilian space," as Blumenthal and Kelly proposed
(see \cite{B1959a}). Regarding the terminology, the
reason why
Barbilian proposed this name is the following. In
Barbilian's
construction from the metrization procedure (see
\cite{B1959a}),
consider $J=\mathbb{R}^2$, and let $A,B \in J$. Let $K
\subset
\mathbb{R}^2$ be the Apollonius circle determined by
the points
$A$ and $B,$ and by the constant $\alpha \neq 1,
\alpha \in
\mathbb{R}^*_+.$ This is known to be the geometric
locus of those
points in the plane such that $\frac{PA}{PB}= \alpha,$
for all $P
\in K.$ Then, it is clear that
$$ \mathcal{M} = \max \frac{PA}{PB}= \min
\frac{PA}{PB} = m = \alpha.$$ Therefore, when we
compute the
distance between $A$ and $B,$  $\ln \frac{
\mathcal{M}}{m}=0,$
while $A \neq B.$ This suggested to
Barbilian that,
in order to obtain a distance instead of a
pseudo-distance, one must
avoid the case where $K$ is precisely the Apollonius
circle.
Therefore, he called the distance obtained from the
above
pseudo-metric {\it Apollonian metric} or {\it Apollonian
distance.}
 This terminology is still used in many references
today.

It is our hope that the chronological details
presented in this
note supplement the information already available in
this area. A historical note on this topic written by the authors of the present material
will appear in {\it Historia Mathematica} \cite{BS2007}. As
for Barbilian's life, literary work or other
biographical details,
there are many references available to interested
readers, as for
example the comprehensive monograph in English \cite{C1981}, as 
well as other sources in Romanian \cite{C1989,G1979}.

The authors express their thanks to Alfonso Agnew
 for his many
useful suggestions that improved this presentation.

 \vspace{.2cm}

{\it Authors' adresses:} \\
W. G. Boskoff, Department of
Mathematics
and Computer
Science, University Ovidius, Constantza, Romania\\
email: boskoff@univ-ovidius.ro

\vspace{.25cm}
\noindent
B. D. Suceav\u{a}, 
Department of Mathematics, California State Unversity
Fullerton,
P.O. Box 6850, Fullerton,
CA 92834-6850\\
email: bsuceava@fullerton.edu

\begin{thebibliography}{B1934}
 



\bibitem{B1934} D. Barbilian - {\it Einordnung von
Lobayschewskys Massenbestimmung in einer gewissen
allgemeinen
Metrik der Jordansche Bereiche,} \v{C}asopsis Mathematiky
a Fysiky
{\bf 64} (1934-35), 182-183.

\bibitem{B1933-1934} D. Barbilian - {\it Curs de
matematici elementare \c{s}i geometrie
descriptiv\u{a}}, Editura
Facult\u{a}\c{t}ii de \c{s}tiin\c{t}e, Bucure\c{s}ti,
1933-1934;
published in {\it Opera didactic\u{a}}, volume I,
editor N.
Mih\u{a}ileanu, foreword by G. Vr\u{a}nceanu.

\bibitem{B1940} D. Barbilian - {\it Not\u{a} asupra
lucr\u{a}rilor \c{s}tiin\c{t}ifice,} Ed. Bucovina,
I.E.
Torou\c{t}iu, 1940.

\bibitem{B1959a} D. Barbilian - {\it Asupra unui
principiu de metrizare,} Stud. Cercet. Mat. {\bf 10}
(1959), 68
- 116.

\bibitem{B1959b} D. Barbilian - {\it Fundamentele
metricilor abstracte ale lui Poincar\'e \c{s}i
Carath\'eodory ca
aplica\c{t}ie a unui principiu general de metrizare,}
Stud.
Cercet. Mat. {\bf 10} (1959), 273 - 306.

\bibitem{B1960} D. Barbilian - {\it J-metricile
naturale finsleriene,} Stud. Cercet. Mat. {\bf 11}
(1960), 7 -
44.

\bibitem{B1968} D. Barbilian - {\it Geometrie - Opera
didactic\u{a},} Ed. Tehnic\u{a}, Bucure\c{s}ti, 1968.

\bibitem{BR1962} D. Barbilian and N. Radu - {\it
J-metricile naturale finsleriene \c{s}i func\c{t}ia de
reprezentare a lui Riemann,} Stud. Cercet. Mat. {\bf
12} (1962),
pp. 21 - 36.

\bibitem{Beardon1998} A. F. Beardon - {\it The
Apollonian metric of a domain in $\mathbb{R}^n$,} in
{\it
Quasiconformal mappings and analysis,}
Springer-Verlag, 1998,
91-108.

\bibitem{B1938} L. M. Blumenthal - {\it Distance
Geometry. A Study of the Development of Abstract
Metrics. With an
introduction by K. Menger}, Univ. of Missouri Studies
13, Univ. of
Missouri, Columbia, 1938.

\bibitem{Blum1954} L. M. Blumenthal - Review of {\it 
Kelly, P. J. 
 Barbilian geometry and the Poincar\'e model},  
Amer. Math. Monthly  61,  (1954). 311-319, Mathematical Reviews / MathSciNet, MR0061397 (15,819a).

\bibitem{Blum1959} L. M. Blumenthal - Review of {\it Barbilian, D. 
 Sur un principe de m\'etrisation}, Mathematical Reviews / MathSciNet,
MR0107848 (21 \# 6570). 

\bibitem{WB1993} W.  G. Boskoff - {\it About metric
Barbilian spaces,} Sci. Bull., Politeh. Univ. Buchar.,
Ser. A {\bf
55} (1993) No.3-4, 61 - 70.

\bibitem{WB1994a} W. G. Boskoff and P. Horja - {\it
$S$-Riemannian manifolds and Barbilian spaces,} Stud.
Cercet. Mat.
{\bf 46} (1994), No.3, pp. 317 - 325.

\bibitem{WB1994b} W. G. Boskoff - {\it The
characterization of some spectral Barbilian spaces
using the
Tzitzeica construction,} Stud. Cercet. Mat. {\bf 46}
(1994) No.5,
pp. 503 - 514.

\bibitem{WB1995} W. G. Boskoff - {\it Finslerian and
induced Riemannian structures for natural Barbilian
spaces,} Stud.
Cercet. Mat. {\bf 47} (1995) No.1, pp. 9-16.

\bibitem{WB1996} W. G. Boskoff - {\it The connection
between Barbilian and Hadamard spaces}, Bull. Math.
Soc. Sci.
Math. Roum., Nouv. SZr. {\bf 39} (1996) No.1-4, pp.
105 - 111.

\bibitem{WB1996carte} W. G. Boskoff - {\it Hyperbolic
geometry and Barbilian spaces,} Istituto per la
Ricerca di Base,
Hardronic Press, 1996.

\bibitem{WB1998} W. G. Boskoff - {\it A generalized
Lagrange space induced by the Barbilian distance,}
Stud. Cercet.
Mat. {\bf 50} (1998) No.3-4, pp. 125-129.

\bibitem{WB2002} W. G. Boskoff - {\it Manifolds with
Barbilian metric structure. (Varietati cu structura
metrica
Barbilian). (Romanian)} Colec\c{t}ia Biblioteca de
Matematic\u{a},
Ed. Ex Ponto, Constanta, 2002.

\bibitem{BS2007} W. G. Boskoff, B. D. Suceav\u{a} - {\it  Barbilian Spaces: The History of a Geometric Idea}, Historia Mathematica (to appear).

\bibitem{BGS2007} W. G. Boskoff, M. G. Ciuc\u{a},
B. D. Suceav\u{a}- {\it Distances Induced by Barbilian's Metrization
Procedure}, Houston J. Math. (to appear). LANL math.DG/0606608.

\bibitem{C1966} M. C\u{a}linescu - {\it Insemn\u{a}ri despre motivul poetic al oglinzii 
(Plec\^and de la
Jocul secund al lui Ion Barbu)}, Via\c{t}a Rom\^aneasc\u{a}, {\bf 11} (1966).

\bibitem{C1981} A. Cior\u{a}nescu - {\it Ion Barbu}, Twayne Publishers, Boston, 1981.

\bibitem{C1989} M. Colo\c{s}enco - {\it Ion Barbu -
Dan Barbilian}, Ed. Minerva, 1989.

\bibitem{D1984} D. Pillat - {\it Foreword} to the edition {\it Ion Barbu - Versuri \c{s}i proz\u{a}},
Editura Minerva, Bucure\c{s}ti, 1984.

\bibitem{GH2000} F. W. Gehring and K. Hag - {\it  The
Apollonian metric and quasiconformal mappings},
Contemp. Math.,
256 (2000) 143--163.

\bibitem{G1979} G. Gibescu - {\it Chronological
Note}, foreword to {\it Ion Barbu - Poezii,} Ed.
Albatros, 
Bucharest, 1979.

\bibitem{H2003} P. A. H\"ast\"o - {\it The
Appollonian metric: uniformity and quasiconvexity},
Ann. Acad.
Sci. Fennicae, {\bf 28} (2003), pp. 385 - 414.

\bibitem{H2003b}  P. A. H\"ast\"o - {\it The
Apollonian metric: limits of the approximation and
bilipschitz
properties,} Abstr. Appl. Anal. 2003, no. 20,
1141-1158.

\bibitem{H2004} P. A. H\"ast\"o - {\it The Apollonian
inner metric,} Comm. Anal. Geom. {\bf 12} (2004),
927-947.

\bibitem{H2004b} P. A. H\"ast\"o - {\it The Apollonian
metric: quasi-isotropy and Seittenranta's metric,}
Comput. Methods
Funct. Theory {\bf 4} (2004), no. 2, 249-273.


\bibitem{H2005b} P. A. H\"ast\"o - {\it  The Apollonian metric: the comparison property, bilipchitz mappings and
 thick sets}, J. Appl. Anal. 12 (2006), no. 2 (to appear).

\bibitem{HI2005} P. A. H\"ast\"o and Z. Ibragimov - {\it Apollonian
isometries of planar domains are M\"{o}bius mappings}, J. Geom. Anal. 
15 (2005), no. 2, 229-237.



\bibitem{HL2004}  P. A. H\"ast\"o and H. Lind\'en - {\it Isometries of the half-apollonian metric}, Complex Var. Theory Appl. {\bf 49}
 (2004), 405-415. 

\bibitem{HPS2005} P. A. H\"ast\"o, S. Ponnusamy and S. K. Sahoo
- {\it  Inequalities and geometry of the Apollonian and related metrics} Rev. Roumaine Math. Pures Appl. (to appear)

\bibitem{I2002} Z. Ibragimov - {\it The Apollonian
metric, sets of constant width and M\"obius modulus of
ring
domains,} Ph.D. Thesis, University of Michigan, Ann
Arbor, 2002.

\bibitem{I2003a} Z. Ibragimov - {\it
On the Apollonian metric of domains in
$\overline{\mathbb{R}}{}\sp
n$,} Complex Var. Theory Appl. {\bf 48} (2003), no.
10, 837-855.

\bibitem{I2003b} Z. Ibragimov - {\it Conformality of
the Apollonian metric,} Comput. Methods Funct. Theory,
 {\bf 3}
(2003), 397-411.

\bibitem{K1954} P. J. Kelly  - {\it Barbilian
Geometry and the Poincar\'e Model}, Amer. Math.
Monthly, {\bf 61}
(1954), pp. 311 - 319.

\bibitem{KM1981} P. J. Kelly and G. Matthews - {\it
The Non-Euclidean, Hyperbolic Plane. Its structure and
consistency} Springer-Verlag, 1981.

\bibitem{K1998} Sh. Kobayashi - {\it Hyperbolic Complex Spaces,} 
Springer-Verlag, 1998.

\bibitem{N1960} B. Nicolescu - {\it Ion Barbu. Cronologia Jocului secund}, Editura pentru Literatur\u{a}, Bucure\c{s}ti, 1960.

\bibitem{S1999} P. Souza -
{\it Barbilian metric spaces and the hyperbolic plane.} (Spanish)
 Misc. Mat. {\bf 29,} (1999), 25-42. 

\bibitem{V1935} T. Vianu - {\it Ion Barbu}, Cultura Na\c{t}ional\u{a}, Bucure\c{s}ti, 1935.

\end{thebibliography}
\end{document}